\documentclass[12pt,a4paper]{amsproc}
\usepackage[english]{babel}
\usepackage{latexsym}
\usepackage{amsmath}
\usepackage{amssymb}
\usepackage[all]{xy}

\usepackage{amsmath,tikz-cd}
\usepackage{amscd}
\usepackage[cp850]{inputenc}
\usepackage[mathscr]{eucal}
\newcommand{\To}{\longrightarrow}
\newcommand{\bil}[2]{{\left\langle\kern0ex #1,#2
		\kern0ex\right\rangle}}

\numberwithin{equation}{section}

\newtheorem{thm}{Theorem}[section]

\newtheorem{lem}[thm]{Lemma}
\newtheorem*{claim}{Claim}
\newtheorem{prop}[thm]{Proposition}
\newtheorem{quest}{Question}[section]

\newtheorem{defin}[quest]{Definition}
\newtheorem{theorem}[quest]{Theorem}

\newtheorem{corollary}[quest]{Corollary}

\theoremstyle{remark}

\newcommand{\N}{\mathbb{N}}

\newcommand{\adef}{\begin{defin}}
\newcommand{\zdef}{\end{defin}}

\def\PB{\operatorname{PB}}
\def\PO{\operatorname{PO}}
\newcommand\restr[2]{{
		\left.\kern-\nulldelimiterspace 
		#1 
		\right|_{#2} 
}}

\newcommand{\aproof}{\begin{proof}}
\newcommand{\zproof}{\end{proof}}

\title{Polyhedrality for twisted sums with $C(\omega^\alpha)$}

\author{Jes\'{u}s M. F. Castillo}
\address{Universidad de Extremadura\\Instituto de Matem\'aticas Imuex\\
Avenida de Elvas s/n\\ 06071-Badajoz\\ Spain} \email{castillo@unex.es}

\author{Alberto Salguero Alarc\'on}
\address{Universidad Complutense de Madrid\\ Plaza de las Ciencias, s/n \\ 28040-Madrid\\ Spain} \email{albsalgu@ucm.es}



\thanks{This research was supported in part by MINCIN project PID2019-103961GB-C21 and PID2023-146505NB-C21, and by Junta de Extremadura project IB20038. The second-named author benefited from an FPU Grant FPU18/00990 from the Spanish Ministry of Universities.}

\thanks{The authors thank the anonymous referee for a very careful reading pointing out a few mistakes and possible improvements in the manuscript. Definitely, this version is better.}

\begin{document}

\maketitle

\begin{abstract} We obtain two partial answers to the $3$-space problem for isomorphic polyhedrality: (1) every twisted sum of $C(\alpha)$, $\alpha<\omega_1$,  with a separable isomorphically polyhedral space with the BAP, is isomorphically polyhedral. (2) Every twisted sum of $c_0(\aleph)$ and a Banach space having a boundary with property $(*)$ has a boundary with property $(*)$, hence it is isomorphically polyhedral.\end{abstract}

\section{Introduction and preliminaries}
\par This paper is concerned with the stability of polyhedral renormings under exact sequences of Banach spaces. Recall that a Banach space is \emph{polyhedral} if every finite dimensional (equivalently, $2$-dimensional) subspace is isometrically isomorphic to a subspace of some finite-dimensional $\ell_\infty^n$. Banach spaces admitting a polyhedral renorming are called \emph{isomorphically polyhedral}.  A survey-like reference about open problems regarding polyhedrality is \cite{troy}. Besides, an \emph{exact sequence} of Banach spaces is a diagram of Banach spaces and (linear, continuous) operators
\begin{equation} \label{eq:ex-seq} \xymatrix{0 \ar[r] & Y \ar[r]^j & Z\ar[r]^q  & X \ar[r] & 0} \end{equation}
 in which the kernel of every operator agrees with the image of the preceding. This amounts to saying that $j$ is an into isomorphism, $q$ is a quotient map and $X$ is isomorphic to $Z/j(Y)$. The middle space $Z$ is called a \emph{twisted sum} of $Y$ and $X$. The by now standard \emph{$3$-space properties} of Banach spaces are those which are stable under twisted sums; that is to say, a property $\mathcal P$ is a 3-space property if whenever the spaces $Y$ and $X$ in an exact sequence of Banach spaces (\ref{eq:ex-seq}) have property $\mathcal P$ then also $Z$ has property $\mathcal P$.
The monograph \cite{castgonz} is a good reference about $3$-space properties.

In many times and places has the first author asked \cite{troy,castpapi} whether to be isomorphically polyhedral is a 3-space property. The paper \cite{csQM} provided the first partial answer showing that every twisted sum of two $c_0(\Gamma)$ spaces is isomorphically polyhedral. The argument there presented contains however an incorrect assertion: we will deal with this in Section \ref{boundaries}, where we additionally present an improved version for that result. After $c_0(\Gamma)$, the next simplest isomorphically polyhedral spaces are the spaces of continuous functions on some countable ordinal space. We provide in Section 2 another positive partial result in that direction, showing that every twisted sum of $C(\alpha)$, where $\alpha<\omega_1$ and an isomorphically polyhedral space $X$ is again isomorphically polyhedral provided $X$ is separable with the BAP.

\subsection{Background on exact sequences}
We provide here a brief background on the theory of exact sequences in order to make this paper reasonably self-contained. We address the interested reader to the monograph \cite{hmbst} for a fully detailed exposition.
\par It is by now a classical result \cite{kalton} that exact sequences can be represented by \emph{quasilinear maps}, which are homogeneous maps $\Omega: X \to Y$ such that $\|\Omega(x+y)-\Omega(x)-\Omega(y)\| \leq M(\|x\|+\|y\|)$.
Indeed, every such $\Omega$ induces an exact sequence
\begin{equation} \label{eq:ex-quasi} \xymatrix{0 \ar[r] & Y \ar[r]^-{\iota} & Y\oplus_\Omega X \ar[r]^-{\pi}  & X \ar[r] & 0} \end{equation}
where $Y\oplus_\Omega X$ is the vector space $Y\times X$ endowed with the quasi-norm $\|(y,x)\|_\Omega = \|y-\Omega(x)\|+\|x\|$, while $\iota(y) = (y,0)$ and $\pi(y,x)=x$. Moreover, every exact sequence (\ref{eq:ex-seq}) is equivalent to an exact sequence (\ref{eq:ex-quasi}), in the sense that there exists an isomorphism $u$ making commutative the following diagram
\[  \xymatrix{0 \ar[r] & Y \ar[r]^j \ar@{=}[d] & Z\ar[r]^q \ar[d]^u & X \ar@{=}[d] \ar[r]& 0 \\ 0 \ar[r] & Y \ar[r]^-{\iota} & Y\oplus_\Omega X \ar[r]^-{\pi}  & X \ar[r] & 0}
\]
An exact sequence (\ref{eq:ex-seq}) is \emph{trivial}, or it \emph{splits}, if $j(Y)$ is complemented in $Z$. In terms of quasilinear maps, this happens if and only if the corresponding map $\Omega$ is \emph{trivial}, with the meaning that it can be written as $\Omega = L + B$, where $L: X\to Y$ is a (not necessarily continuous) linear map and $B: X\to Y$ is a bounded homogeneous map (there is $C>0$ such that $\|Bx\|\leq C\|x\|$. The infimum of those $C$ will be denoted $\|B\|$). We will write $\Omega \equiv 0$ to denote that $\Omega$ is trivial.
At some point it will be important to quantify triviality: If $\Omega = L + B$ then we will say that $\Omega$ is $\|B\|$-trivial. A sequence $\Omega_n$ of trivial maps will be called \emph{uniformly trivial} if the decompositions $\Omega_n = B_n + L_n$ can be chosen such that $\sup_n \|B_n\|<\infty$.
We will say that $\Omega$ and $\Phi$ are equivalent, and write $\Omega\equiv \Phi$,  to mean $\Omega - \Phi \equiv 0$, and so two sequences of quasilinear maps $\Omega_n$ and $\Phi_n$ are \emph{uniformly equivalent} when $\Omega_n - \Phi_n$ is uniformly trivial. A crucial point in the theory of twisted sums is that the space $Z$ does not need to be locally convex; it will actually be isomorphic to a Banach space precisely when $\Omega$ is \emph{1-quasilinear}, which means that there is $M>0$ such that for every $n\in \N$ and $x_1, ..., x_n\in X$ one has $\|\Omega(\sum_{i=1}^n x_i)-\sum_{i=1}^n\Omega(x_i)\| \leq M(\sum_{i=1}^n\|x_i\|)$.

\par Finally, given an exact sequence (\ref{eq:ex-seq}) and an operator $T: W \to X$, we can produce the lower exact sequence in the diagram
\[\xymatrix{0 \ar[r] & Y \ar[r]^j & Z \ar[r]^q & X \ar[r] & 0 \\
0 \ar[r] & Y \ar[r]^{\underline \jmath} \ar@{=}[u] & \PB \ar[r]^{\bar q} \ar[u]^{\underline T} & W \ar[r] \ar[u]^T & 0 }
\]

\noindent called the \emph{pullback sequence}, in which $PB = \{(w,z)\in W \oplus_\infty Z: T(w) = q(z)\}$. The operators $\underline q$ and $\underline T$ are just the restrictions to $PB$ of the canonical projections of $W\oplus_\infty Z$ to $W$ and $Z$, respectively, and $\underline\jmath(y) = (j(y),0)$. It is straightforward to see that if $T$ is an into isomorphism, then so is $\underline T$. The pullback sequence is trivial precisely when $T$ admits a \emph{lifting} through $q$; that is, an operator $\hat{T}: W \to Z$ such that $q\hat{T}=T$.
\par Analogously, given an exact sequence (\ref{eq:ex-seq}) and an operator $T: Y \to V$, the lower exact sequence in the diagram
\[\xymatrix{0 \ar[r] & Y \ar[r]^j  \ar[d]^{T} & Z \ar[r]^q \ar[d]^{\overline T} & X \ar[r] \ar@{=}[d] & 0 \\
	0 \ar[r] & V \ar[r]^{\overline\jmath} & \PO \ar[r]^{\overline q}  & X \ar[r] & 0 }
\]

\noindent is known as the \emph{pushout sequence}. The space $PO$ is the quotient of $Z\oplus_1 V$ by the closure of $\Delta = \{(jy, -Ty): y \in Y)\}$, the operators $\overline{T}$ and $\overline{\jmath}$ are, respectively, the canonical inclusions of $Z$ and $V$ into $Z\oplus_1 V$ composed with the quotient map $Z\oplus_1 V \to PO$, and $\overline{q}[(z,v)+\overline{\Delta}] = q(z)$. Again, if $T$ is an into isomorphism, then so is $\overline{T}$. Finally, the pushout sequence splits if and only if $T$ admits an \emph{extension} to $Z$, with the meaning that there is $\tilde{T}: Z \to V$ such that $\tilde Tj = T$.

\section{The $3$-space problem for isomorphic polyhedrality}

In this section we prove one of the main results in the paper. In its general form we get:

\begin{prop} Let $X, Y_N$ be separable isomorphically polyhedral spaces, $X$ with the BAP, and let $Y$ be a subspace of $c_0(\mathbb N, Y_N)$. If every twisted sum of $Y_N$ and $X$ is isomorphically polyhedral then every twisted sum of $Y$ and $X$ is isomorphially polyhedral.\end{prop}
\begin{proof} Let $Z$ be a twisted sum as indicated, namely, the middle space in an exact sequence
\begin{equation}\label{seq}\xymatrix{0 \ar[r] & c_0(\mathbb N, Y_N) \ar[r] & Z\ar[r]  & X \ar[r] & 0}\end{equation}
Let us call $\Omega$ the $1$-quasilinear map associated to the sequence (\ref{seq}), which means that
(\ref{seq}) is equivalent to the exact sequence
$$\xymatrix{0 \ar[r] & c_0(\mathbb N, Y_N) \ar[r] &  c_0(\mathbb N, Y_n)\oplus_\Omega X \ar[r]  & X \ar[r] & 0}.$$

Let us consider $P_N: c_0(\mathbb N, Y_N) \To Y_N $ the canonical contractive projection. We invoke now the ``chasing device" of \cite{castmoresob} (see also \cite[Proposition 10.1.3]{hmbst}) to obtain that, in the presence of the BAP, if we call $\Phi:  X\To c_0(\mathbb N, Y_N))$ the $1$-quasilinear map $$\Phi (x) = (P_N \Omega x)_{N=1}^\infty$$
then $\Omega -  \Phi$ is trivial because its components $P_N\Omega - P_N\Omega =0$ are uniformly trivial. Consider the exact sequences
$$\xymatrix{0 \ar[r] &  Y_N  \ar[r]^-{\imath_N} & Y_N \oplus_{P_N\Omega} X \ar[r]^-{Q_N}  & X \ar[r] & 0}$$
make their $c_0$-sum \cite[Prop. 10.2.1]{hmbst}
\begin{equation}\label{amalgam} \xymatrix{0 \ar[r] &  c_0(\mathbb N, Y_N) \ar[r]^-{(\imath_N)} &  c_0\Big(\N, Y_N \oplus_{P_N\Omega} X \Big)\ar[r]^-{(Q_N)}  & c_0(\N, X) \ar[r] & 0}\end{equation}
and call $c_0(P_N\Omega)$ the $1$-quasilinear map that defines the exact sequence (\ref{amalgam}). Now, look in awe at the diagram

$$\xymatrix{0 \ar[r] &  c_0(\mathbb N, Y_N) \ar@{=}[d]  \ar[r]^-{(\imath_N)} &  c_0(\mathbb N, Y_N)  \oplus_{P_N\Omega} X )\ar[r]  & c_0(\N, X) \ar[r] & 0\\
0 \ar[r] &  c_0(\mathbb N, Y_N)  \ar[r] & Z \ar[r] & X \ar[r] & 0}
$$
and recall that the $1$-quasilinear map defining the lower sequence is $\Phi$.

\begin{claim} $ (\imath_N)\Phi\equiv 0$.
\end{claim}

\begin{proof}[Proof of the Claim] Indeed, $(\imath_N)\Phi$ generates an exact sequence
\[\xymatrix{0 \ar[r] &  c_0\big (\mathbb N, Y_N \oplus_{P_N\Omega} X \Big)\ar[r]  & M \ar[r] & X \ar[r]& 0}.
\]
Since $\imath_N \Phi = \imath_N P_N\Omega$, it is clear that the components of $(\imath_N)\Phi$ are uniformly trivial. More precisely, if $\Psi$ is a quasilinear map, and we consider the exact sequence $\xymatrix{0 \ar[r] &  U\ar[r]^-\imath & U\oplus_{\Psi} W \ar[r]  & W\ar[r]& 0}$
then $\imath \Psi $ is trivial since the linear map $L: W \to U\oplus_{\Psi} W$ given by $L(w)= - (0,w)$ is at distance $1$ from $\imath \Psi$:
$$\|\imath \Psi (w)  - L(w)\| = \|(\Psi w, 0) + (0,w)\| = \|(\Psi w, w)\| =\|w\|$$

A second invocation to the chasing device yields the Claim.\end{proof}

The inexorable consequence of the Claim is that by \cite[Proposition 2.11.3]{hmbst} the lower sequence is a pullback of the upper sequence, which means the existence of a commutative  diagram

$$\xymatrix{0 \ar[r] &  c_0(\mathbb N, Y_N) \ar@{=}[d]  \ar[r]^-{(\imath_N)} &  c_0\Big (\mathbb N, Y_N \oplus_{P_N\Omega} X \Big )\ar[r]^-{(Q_N)}  & c_0(\N, X) \ar[r] & 0\\
0 \ar[r] &  c_0(\mathbb N, Y_N) \ar[r] & Z \ar[u]_v \ar[r] & X \ar[r] \ar[u]_u& 0}
$$
Therefore, $Z$ is isomorphic to the pullback space, namely, the subspace $\{( (f_N, x):\; (Q_N f_N)= u x\}$ of the product $ c_0\Big (\mathbb N, Y_N \oplus_{P_N\Omega} X \Big )\oplus_\infty X$. Now, the space $Y_N \oplus_{P_N\Omega} X$ is isomorphically polyhedral by hypothesis, and the main result in \cite{dfh} asserts that in a separable isomorphically polyhedral space every norm can, for every $\varepsilon >0$, be $(1+\varepsilon)$-approximated by a polyhedral norm, say $\|\cdot\|_N$. Therefore, we could well assume that $Z$ is a subspace of $c_0 (\N, Y_N) \oplus_{\|\cdot\|_N} X)\oplus_\infty X$. Since any $c_0$-sum of polyhedral spaces is polyhedral \cite{hn}, the space $Z$ would be isomorphically polyhedral.\medskip

The assertion about the subspace $Y$ follows from the previous theorem, the pushout diagram
	\[\xymatrix{0 \ar[r] & Y \ar[r]  \ar[d]^{\iota} & Z \ar[r] \ar[d]^{\overline\iota} & X \ar[r] \ar@{=}[d] & 0 \\
		0 \ar[r] & C(\alpha) \ar[r] & PO \ar[r] & X \ar[r] & 0 }
	\]
	and the fact that if $\iota$ is an into isomorphism, then so is $\overline\iota$.\end{proof}

The specific form we aim to is:

\begin{theorem} Let $\alpha<\omega_1$. If $X$ is an isomorphically polyhedral separable space with the BAP then every
twisted sum of $C(\alpha)$ and $X$ is isomorphically polyhedral.\end{theorem}
\begin{proof} We just have to check that we are in the hypotheses of the Proposition above. Since $C(\alpha)\simeq C(\alpha +1)$ it is clear that only the limit ordinal case needs to be considered, and when $\alpha=\lim \alpha_N$ is a limit ordinal $C(\alpha) \simeq c_0(\N, C(\alpha_N)$.
To check that every twisted sum of $C(\alpha_N)$ and $X$ is isomorphically polyhedral we will just proceed by transfinite induction observing that for $\alpha_N = \omega^N$ the space $C(\omega^N)\simeq c_0$ is $(2N+1)$-\emph{separably injective} \cite[Prop. 2.23]{2132}, every twisted sum of $C(\omega^N)$ and $X$ is $(2N+1)$-isomorphic to $C(\omega^N)\oplus_\infty X$ and therefore it is isomorphically polyhedral.\end{proof}

\section{Boundaries}\label{boundaries}

Let us report an incorrect assertion in \cite{csQM} and amend and improve the statement of \cite[Theorem 4.2]{csQM}.

\par The claim to be corrected, \cite[p.13]{csQM}, says that a Banach space with a $\sigma$-discrete boundary is isomorphically polyhedral. Recall that given $X$ a Banach space, a subset $B$ of $S_{X^*}$ is a \emph{boundary} if for every $x\in X$ there is $x^*\in B$ so that $x^*(x) = \|x\|$. A boundary is said to be \emph{$\sigma$-discrete} if it can be written as a countable union of relatively discrete sets. If \emph{discrete} here means ``discrete in the norm topology'', then the assertion in \cite[p.13]{csQM} is wrong, and $\ell_1$ is the perfect example: this space is not isomorphically polyhedral since it does not contain $c_0$ while $\{+1,-1\}^\N$ is a norm-discrete weak*-compact boundary (hence it cannot be discrete in the weak*-topology). It is clear that a Banach space that admits a weak*-$\sigma$-discrete boundary also admits a $\sigma$-discrete in the norm topology boundary, although the converse fails: to prove it, recall that a point $f\in S_{X^*}$ is called a \emph{$w^*$-exposed point}  if there exist $p_f\in S_X$ such that $f(p_f)=1$ and $g(p_f)<1$ for every $g\in S_{X^*}$ different from $f$. It is clear that the set $\checkmark$ of $w^*$-exposed points is contained in any boundary and also in the set of extreme points.

\begin{lem} $\ell_1$ does not admit any  weak*-$\sigma$-discrete boundary.
\end{lem}

\begin{proof} It suffices to recall that, in this case, $\checkmark$ coincides with the set of extreme points $\{+1,-1\}^\N$ \cite[Example 1.15]{bible}.\end{proof}

Having a norm discrete boundary $\{+1,-1\}^\N$ but not a weak*-$\sigma$-discrete boundary does not happen in isomorphically polyhedral spaces. To prove it we recall another boundary: the set $\S$ of $w^*$-strongly exposed points, where a point $f\in S_{X^*}$ is called a \emph{$w^*$-strongly exposed point} if there exist $p_f\in S_X$ such that $f(p_f)=1$ and whenever $(g_n)\subset B_{X^*}$ is a sequence such that $\lim g_n(p_f)=1$ then $\lim \|g_n - f\|=0$. Observe that $\S\subset \checkmark$. We have:

\begin{lem} An isomorphically polyhedral Banach space admits a  weak*-$\sigma$-discrete boundary  if and only if it admits a  norm-$\sigma$-discrete boundary  \end{lem}
\begin{proof} A beautiful result of Fonf \cite{fonfext} (see also \cite{vesselyext}) establishes that if $X$ is a polyhedral space then
$\S$ is a boundary (hence a minimal boundary). Thus, if $X$ admits a boundary that is discrete in the weak*-topology, then $\S$ must be
discrete in the norm-topology as well. Therefore, given $f\in \S$, no sequence $(g_n)$ of $\S$ is such that $\lim \|g_n -f\|=0$; consequently $\lim g_n(p_f)=1$ is impossible too, which means that $f$ is weak*-isolated in $\S$.\end{proof}

In \cite[Theorem 4.2]{csQM}, the authors showed that every twisted sum of $c_0(\kappa)$ and a Banach space $X$ having a weak$^*$-$\sigma$-discrete boundary has itself a weak$^*$-$\sigma$-discrete boundary. This part of the proof is correct. However, the coda in the statement \cite[Theorem 4.2]{csQM} \emph{In particular, every twisted sum of $c_0(\kappa)$ and a Banach space having a weak$^*$-$\sigma$-discrete boundary is isomorphically polyhedral} is uncertain because we do not know if the existence of a boundary that is weak*-$\sigma$-discrete implies that the space is isomorphically polyhedral. A closely related result is however possible. Recall that a boundary $\mathfrak B$ is said having \emph{property (*)} if whenever $f$ is a weak*-accumulation point of $\mathfrak B$ and $x_0$ is a normalized element then $f(x_0)<1$. Banach spaces with boundary enjoying property $(*)$ form the most general class known nowadays of isomorphically polyhedral spaces. In fact, a separable Banach space is isomorphically polyhedral if and only if it admits a boundary having property $(*)$ \cite{fonf-russian}. A verbatim following of \cite[Theorem 11]{fpst} yields:

\begin{prop} If a space $X$ possesses a boundary which is \emph{weak*-compact} and weak*-$\sigma$-discrete then it admits for every $\varepsilon>0$ an $\varepsilon$-equivalent polyhedral renorming having a boundary with property $(*)$.\end{prop}

\par After all those prolegomena, let us show that the result announced in \cite[Corollary 4.3]{csQM}, namely, that every twisted sum of $c_0(\kappa)$ and $c_0(I)$ is isomorphically polyhedral, remains true. One actually has:

\begin{prop} If $X$ is a Banach space having a boundary $\mathfrak B$ with property $(*)$ then every twisted sum of $c_0(\kappa)$ and $X$ has a boundary with property $(*)$, hence it is isomorphically polyhedral.\end{prop}
\begin{proof}  The key tool is the representation theorem \cite[Theorem 3.1]{csQM}. It asserts that for every exact sequence $\xymatrix{0 \ar[r] & c_0(\kappa) \ar[r]^j & Z\ar[r]^q  & X \ar[r] & 0}$ there is a commutative diagram
	\begin{equation}\label{diagram} \xymatrix{0 \ar[r] & c_0(\kappa) \ar[r] & C(L) \ar[r]^r& C(B_{X^*}) \ar[r]  & 0 & \\
0 \ar[r] & c_0(\kappa) \ar[r]^j \ar@{=}[u] & Z\ar[r]^q \ar[u] & X \ar[r] \ar[u]_{\delta_X} & 0}
			\end{equation}
where $L$ is a \emph{$\kappa$-discrete extension} of $(B_{X^*}, w^*)$, that is, $L$ is a compact space containing a homeo\-morphic to $(B_{X^*}, w^*)$ such that $L\setminus B_{X^*}$ is a discrete space of size $\kappa$. Here $r$ is the natural restriction map, $\delta_X: X \to C(B_{X^*})$ is the canonical evaluation map $\delta_X(x)(x^*)=x^*(x)$ and we therefore can assume that $L= B_{X^*}\cup \kappa$. The magic of the result is that $Z$ can then be, after renorming, identified with the subspace $\{f\in C(B_{X^*}\cup \kappa): \exists x\in X: \; f|_{B_{X^*}} = \delta_X(x)\}$ and therefore the unlabeled arrow is plain inclusion. Let us finally recall that  it was observed in \cite{csQM} that if $\mathfrak B\subset B_{X^*}$ is a boundary for $X$ then $V = \{\delta_\alpha|_{Z} : \alpha <\kappa\} \cup \{\delta_b|_{Z}: b\in \mathfrak B\}$ is a boundary for $Z$.

\par The twist comes now: the set
$$W = \{\delta_\alpha|_Z: \alpha <\kappa\}\cup \{2\delta_b|_Z: b\in \mathfrak B \} $$
is a boundary for the new norm $\nparallel z \nparallel  = \sup \{|\langle z, v\rangle|: \; v\in W \}$ on $Z$. Let us show that this boundary has property $(*)$:

Let $\mu$ be a weak*-accumulation point of $W$. This $\mu$ is either of the form $\delta_{x^*}|_Z$ or of the form $2\delta_{x^*}|_Z$ for a suitable $x^*\in B_{X^*}$.
Pick now $z\in Z$ with $\nparallel  z \nparallel =1$. Since $rz=\delta Q z$, there is $x_z$ such that $z(x^*)=x^*(x_z)$ for all $x^*\in B_{X^*}$. Since $z(b)\leq \frac{1}{2}$ for all $b\in \mathfrak B$ it follows that $z(x^*)\leq \frac{1}{2}$ for all $x^*\in B_{X^*}$. Now, if $\mu= \delta_{x^*}|_Z$ then $\langle \mu, z\rangle = z(x^*) \leq 1/2<1$. If, however, $\mu = 2\delta_{x^*}|_Z$ this can only happen because
$x^*$ is a weak*-accumulation point of $\mathfrak B$, in which case $\langle \mu, z\rangle = z(x^*) = x^*(x_z)<1$ since $\mathfrak B$ has property $(*)$. \end{proof}


\begin{thebibliography}{123}

\bibitem{2132} A.~Avil\'es, F.~Cabello  S\'anchez, J.M.F.~Castillo, M.~Gonz\'alez,
Y.~Moreno, \emph{Separably injective Banach spaces}, Lecture Notes in Math. 2132, Springer, 2016.


\bibitem{bible} V. Bible, \emph{Discreteness properties, smoothness and polyhedrality}, Doctoral Thesis, University College, Dublin 2016.

\bibitem{hmbst} F. Cabello S\'anchez, J.M.F. Castillo, \emph{Homological methods in Banach space theory}, Cambridge Studies in Advanced Math. 203, (2023) Cambridge Univ. Press. ISBN. 9781108778312.

\bibitem{ccky} F. Cabello S\'{a}nchez, J.M.F. Castillo, N.J. Kalton, D.T. Yost, \emph{Twisted sums with $C(K)$ spaces}, Trans. Amer.
Math. Soc. 355 (2003) 4523--4541.



\bibitem{castgonz} J.M.F. Castillo, M. Gonz\'alez, \emph{Three-space problems in Banach space theory}, Lecture
Notes in Math. 1667, Springer, 1997.
\bibitem{castmoresob} J.M.F.~Castillo, Y.~Moreno, \emph{Sobczyk's theorem and the
Bounded Approximation Property in Banach spaces}, Studia Math. 201
(2010) 1--19.

\bibitem{troy} J.M.F. Castillo, P.L. Papini, \emph{ Hepheastus account on Trojanski's polyhedral war}, Extracta Math. Vol. 29, (2014) 35 -- 51.

\bibitem{castpapi} J.M.F. Castillo, P.L. Papini, \emph{On isomorphically polyhedral $\mathscr L_\infty$-spaces},  J. Funct. Anal. 270 (2016) 2336--2342.

\bibitem{csQM} J.M.F. Castillo, A. Salguero-Alarcón, \emph{Twisted sums of $c_0(I)$}, Quaest. Math. 46 (2023), No. 11, 2339--2354.


\bibitem{dfh} R. Deville, V. Fonf and P. Hajek, \emph{Analytic and polyhedral approximation of convex bodies in
separable polyhedral Banach spaces}, Israel J. Math. 105 (1998)
139 -- 154.

\bibitem{fonf-russian} V.P. Fonf, \emph{Some properties of polyhedral Banach spaces}, Funct. Anal. Appl. 14 (1980) 89--90.

\bibitem{fonfext} V.P. Fonf, \emph{On the boundary of a polyhedral Banach space}, Extracta Math. 15 (2000) 145--154.

\bibitem{fpst}  V.P. Fonf, A.J. Pallares, R.J. Smith, S. Troyanski, \emph{ Polyhedral norms
on non-separable Banach spaces}, J. Funct. Anal. 255 (2008) 449--470.



\bibitem{hn} A.B. Hansen and  N.J. Nielsen, \emph{ On isomorphic classification of polyhedral
preduals  of $L_1$}, Preprint Series Aarhus University, 1973/74
No. 24.

\bibitem{kalton} N.J. Kalton, \emph{The three-space problem for locally bounded F-spaces}, Compo. Math., 37 (1978) 243--276.


\bibitem{vesselyext} L. Vesel\'y, \emph{Boundary of polyhedral spaces: An alternative proof}, Extracta Math. 15 (2000) 213--218.

\end{thebibliography}
\end{document}